\documentclass[11pt, reqno]{amsart}
\usepackage{bbm}
\usepackage{amsfonts}
\usepackage{mathrsfs}
\usepackage{amssymb,amsmath,amscd}
\begin{document}
\renewcommand{\a}{\alpha}
\renewcommand{\b}{\beta}
\renewcommand{\d}{\delta}
\newcommand{\D}{\Delta}
\newcommand{\f}{\frac}
\newcommand{\g}{\gamma}
\newcommand{\G}{\Gamma}
\renewcommand{\l}{\lambda}
\renewcommand{\L}{\Lambda}
\newcommand{\be}{\begin{equation}}
\newcommand{\ee}{\end{equation}}
\newcommand{\bea}{\begin{eqnarray}}
\newcommand{\eea}{\end{eqnarray}}
\newcommand{\bna}{\begin{eqnarray*}}
\newcommand{\ena}{\end{eqnarray*}}
\renewcommand{\o}{\omega}
\renewcommand{\O}{\Omega}
\newcommand{\ov}{\over}
\renewcommand{\le}{\left}
\newcommand{\ri}{\right}
\newcommand{\s}{\sigma}
\renewcommand{\th}{\theta}
\newcommand{\ve}{\varepsilon}
\newcommand{\vp}{\varphi}
\renewcommand{\theequation}{\arabic{section}.\arabic{equation}}
\centerline{\bf\Large Resonance and rapid decay of exponential sums of  }
\centerline{\bf\Large Fourier coefficients of a Maass form for $GL_m(\mathbb Z)$}

\bigskip
\centerline{\sc Xiumin Ren${}^1$ {\it and} Yangbo Ye${}^{1,2}$}

\centerline{${}^1$ Department of Mathematics, Shandong University, Jinan,
Shandong 250100, China}

\centerline{${}^2$ Department of Mathematics, The University of Iowa,
Iowa City, Iowa 52242-1419, USA}

\centerline{Email: {\tt xmren@sdu.edu.cn},
{\tt yangbo-ye@uiowa.edu}}

\bigskip
\centerline{Accepted to appear in SCIENCE CHINA Mathematics}

\bigskip

\date{}
\openup 1\jot
\bigskip

{\footnotesize {\bf  Abstract} \,\,Let $f$ be a
full-level cusp form
for $GL_m(\mathbb Z)$ with Fourier coefficients
$A_f(n_1,\ldots,n_{m-1} )$. In this paper an asymptotic
expansion of Voronoi's summation formula for
 $A_f(n_1,\ldots,n_{m-1} )$ is established. As applications of this
formula, a smoothly weighted average of
$A_f(n,1,\ldots,1)$ against $e(\a |n|^{\b})$ is proved to be rapidly decayed when $0<\b<1/m$. When
$\b=1/m$ and $\a$ equals or approaches $\pm mq^{1/m}$
for a positive integer $q$, this smooth average has a main term of the size of $|A_f(1,\cdots,1, q)+A_f(1,\cdots,1, -q)|X^{1/(2m)+1/2}$,
which is a manifestation of resonance of oscillation exhibited by the Fourier coefficients $A_f(n,1,\cdots,1)$. Similar estimate is also proved for a  sharp-cut sum.}

\bigskip


\bigskip


\medskip
\medskip

\section{Introduction}
\setcounter{equation}{0}

\medskip

Voronoi's summation formula is an important technique in analytic number theory. Ivi\'{c} \cite{Ivi} generated the original Voronoi formula to
non-cuspidal forms as given by multiple divisor functions.
For cuspidal representations of $GL_2(\mathbb Z)$, a Voronoi-type
summation formula was proved by
Sarnak \cite{Sar} for holomorphic cusp forms and by
Kowalski-Michel-Vanderkam \cite{KowMicVan} for Maass forms.
Voronoi's summation formula for Maass forms for
$SL_3(\mathbb Z)$ was proved by Miller-Schmid \cite{MilSch}
and Goldfeld-Li \cite{GdfLi1}. For $m\geq 4$, the Voronoi summation formula was
proved by Miller-Schmid \cite{MilSch1}. To state the formula, let $f$ be a full-level cusp form for $GL_m(\mathbb Z)$  with Langlands'
parameters $\mu_f(j)$, $j=1,\ldots,m$, and Fourier coefficients
$A_f(c_{m-2},\ldots,c_1,n)$. Let
$\psi\in C_c^\infty(\mathbb R^+)$, $q$ a positive integer, and
$h$ an integer coprime with $q$. Let $h\bar{h}\equiv 1 (\bmod\, q)$. Then Voronoi summation formula as
proved by Miller-Schmid \cite{MilSch1} is the following

\begin{eqnarray}
&&
\sum_{n\not=0}
A_f(c_{m-2},c_{m-3},\ldots,c_1,n)
e\Big(-\frac{nh}{q}\Big)
\psi \big(|n|\big)
\nonumber
\\
&=&
q\sum_{d_1|c_1q}
\sum_{d_2|\frac{c_1c_2q}{d_1}}\cdots
\sum_{d_{m-2}|\frac{c_1\cdots c_{m-2}q}{d_1d_2\cdots d_{m-3}}}
\sum_{n\not=0}
\frac{A_f(n,d_{m-2},\ldots, d_1)}{d_1\cdots d_{m-2}|n|}
\nonumber
\\
&\times&
S(n,\bar{h}; q, c, d)
\Psi\Big(\frac{|n|}{q^m}
\prod_{i=1}^{m-2}\frac{d_i^{m-i}}{c_i^{m-i-1}}\Big),
\end{eqnarray}
where $c=(c_1,\ldots, c_{m-2})$, $d=(d_1,\ldots, d_{m-2})$, and
\begin{eqnarray}
S(n, \bar{h}; q, c, d)
&=&
{\mathop{\sum\nolimits^*}
\limits_{x_1\big(\bmod \frac{c_1q}{d_1}\big)}}
e\Big(\frac{d_1x_1n}{q}\Big)
{\mathop{\sum\nolimits^*}
\limits_{x_2\big(\bmod \frac{c_1c_2q}{d_1d_2}\big)}}
e\Big(\frac{d_2x_2\bar{x}_1}{\frac{c_1q}{d_1}}\Big)
\cdots
\nonumber
\\
&\times&
{\mathop{\sum\nolimits^*}
\limits_{x_{m-2}\big(\bmod \frac{c_1\cdots c_{m-2}q}
{d_1\cdots d_{m-2}}\big)}}
e\Big(
\frac{d_{m-2}x_{m-2}\bar{x}_{m-3}}
{\frac{c_1\cdots c_{m-3}q}{d_1\cdots d_{m-3}}}
+
\frac{\bar{h}\bar{x}_{m-2}}
{\frac{c_1\cdots c_{m-2}q}{d_1\cdots d_{m-2}}}
\Big).
\end{eqnarray}
The $*$ in $\sum^*_{t(\bmod r)}$ indicates that $(t,r)=1$.
Here $\Psi$ is an integral transform of $\psi$ given by
\begin{eqnarray}
\Psi(x)
=
\frac1{2\pi i}
\int_{\Re s=-\sigma}\tilde{\psi}(s)x^s
\frac{\tilde{F}(1-s)}{F(s)}ds,
\end{eqnarray}
where
$$
\tilde{\psi}(s)=\int_0^\infty\psi(x)x^s\frac{dx}{x}
$$
and
\begin{eqnarray}
F(s)
&=&
\pi^{-{ms}/2}
\prod_{i=1}^m
\Gamma\Big(\frac{s-\mu_f(j)}{2}\Big),
\\
\tilde{F}(s)
&=&
\pi^{-{ms}/2}
\prod_{i=1}^m
\Gamma\Big(\frac{s-\bar{{\mu}}_f(j)}{2}\Big)
\end{eqnarray}
with
$\{\overline\mu_f(j)\}_{j=1,\ldots,m}
=\{{\mu}_{\tilde{f}}(j)\}_{j=1,\ldots,m}$
being the Langlands'
parameters for the dual form $\tilde{f}$ of $f$.

A special case of (1.1) for even Maass forms for
$SL_m(\mathbb Z)$ was proved by Goldfeld-Li
 \cite{GdfLi1, GdfLi2,GdfLi3}:
\begin{eqnarray}
&&
\sum_{n\not=0}
A_f(1,1,\ldots,1,n)
e\Big(\frac{nh}{q}\Big)
\psi\big(|n|\big)
\nonumber
\\
&=&
q\sum_{d_1|q}
\sum_{d_2|\frac{q}{d_1}}\cdots
\sum_{d_{m-2}|\frac{q}{d_1d_2\cdots d_{m-3}}}
\sum_{n\not=0}
\frac{A_f(n,d_{m-2},\ldots, d_1)}{d_1\cdots d_{m-2}|n|}
\nonumber\\
&\times&
KL(\bar{h},n;d,q)
\Psi\Big(\frac{|n|}{q^m}\prod_{i=1}^{m-2}d_i^{m-i}\Big),
\end{eqnarray}
where $KL(\bar{h},n;d,q)$ is the Kloosterman sum
\begin{eqnarray}
KL(\bar{h},n;d,q)
&=&
{\mathop{\sum\nolimits^*}
\limits_{t_1\big(\bmod \frac{q}{d_1}\big)}}
e\Big(\frac{\bar{h}t_1}{\frac{q}{d_1}}\Big)
{\mathop{\sum\nolimits^*}
\limits_{t_2\big(\bmod \frac{q}{d_1d_2}\big)}}
e\Big(\frac{\bar{t}_1t_2}{\frac{q}{d_1d_2}}\Big)
\cdots
\nonumber\\
&\times&
{\mathop{\sum\nolimits^*}
\limits_{t_{m-2}\big(\bmod \frac{q}{d_1\cdots d_{m-2}}\big)}}
e\Big(\frac{\bar{t}_{m-3}t_{m-2}}{\frac{q}{d_1\cdots d_{m-2}}}\Big)
e\Big(\frac{n\bar{t}_{m-2}}{\frac{q}{d_1\cdots d_{m-2}}}\Big).
\end{eqnarray}
Note that (1.7) can be rewritten as (1.2). When $q=1$, the formulas (1.1) and (1.6) become
\begin{eqnarray}
\sum_{n\neq0}
A_f(1,\ldots,1,n)
\psi\big(|n|\big)
=
\sum_{n\neq0}
\frac{A_f(n,1,\ldots, 1)}{|n|}
\Psi\big(|n|\big).
\end{eqnarray}
Replacing $f$ by its dual form $\tilde{f}$ and noting that
$A_f(n,1,\ldots,1)=A_{\tilde{f}}(1,\ldots,1,n)$, we have
\begin{eqnarray}
\sum_{n\neq0}
A_f(n,1,\ldots,1)
\psi\big(|n|\big)
=
\sum_{n\neq0}
\frac{A_f(1,\ldots,1, n)}{|n|}
\tilde{\Psi}\big(|n|\big),
\end{eqnarray}
where $\tilde{\Psi}(x)$ is defined as $\Psi(x)$ in (1.3) by replacing $\mu_f(j)$ by ${\mu}_{\tilde f}(j)$ in (1.4) and (1.5).
\medskip

In applications, asymptotic behavior of $\Psi(x)$ is
often required.
An asymptotic expansion for Voronoi's summation formula were
firstly obtained by Ivi\'{c} \cite{Ivi} for multiple divisor
functions. Similar asymptotic formulas were proved and used in
subconvexity bounds for $GL_2$ $L$-functions by Sarnak
\cite{Sar} for holomorphic cusp forms and by Liu-Ye
\cite{LiuYe1} and Lau-Liu-Ye \cite{LauLiuYe} for Maass
forms. For Maass forms for $SL_3(\mathbb Z)$, an
asymptotic Voronoi formula was proved by Li \cite{Li} and
Ren-Ye \cite{RenYe4} and applied to subconvexity problems for
$L$-functions attached to a self-dual Maass form for
$SL_3(\mathbb Z)$ by Li \cite{Lii}.

In this paper we prove the following
asymptotic expansion.

\noindent
{\bf Theorem 1.1.} {\it Let $f$ be a full-level cusp form for
$GL_m(\mathbb Z)$. Let $m\geq 3$  be an integer. Let $\Psi(x)$ be as defined in $(1.3)$
with $\psi(y)=\phi(y/X)$, where $\phi(x)\ll 1$ is a fixed smooth function of
compact support on $[a, b]$ with $b>a>0$. Then for any $x>0$, $xX\gg 1$ and $r>m/2$, we have
\begin{eqnarray}
\Psi(x)
&=&
x\sum_{k=0}^{r}
c_k
\int_0^\infty
(xy)^{ 1/(2m)-1/2-k/m}
\psi(y)
\nonumber
\\
&\times&
\Big\{i^{k+(m-1)/2}e\big(m(xy)^{1/m}\big)
+(-i)^{k+(m-1)/2}e\big(-m(xy)^{1/m}\big)\Big\}
dy
\nonumber
\\
&+&
O\big((xX)^{-r/m+1/2+\varepsilon}\big),
\end{eqnarray}
where $c_k$, $k=0,\ldots, r$, are constants depending on $m$ and $\{\mu_f(j)\}_{j=1,\cdots,m}$ with
$c_0=-1/\sqrt{m}$, and the implied constant depends at most on
$f$, $\phi$, $r$, $a$, $b$ and $\varepsilon$.}

We point out that if we replace $f$ be $\tilde{f}$ in (1.3), the formula (1.10) for $\tilde{\Psi}(x)$ has the same leading term, i.e., $\tilde{c}_0=c_0=-1/\sqrt{m}$.
\medskip

As applications of Theorem 1.1,  we consider sums of the Fourier coefficients
$A_f(n,1,\cdots,1)$  twisted with the exponential function $e(\pm\alpha |n|^\beta)$.
Our results are the following.
\medskip

\noindent
{\bf Theorem 1.2.} {\it Let $f$ be a full-level cusp form for
$GL_m(\mathbb Z)$ and $\phi(x)$ a $C^\infty$ function on $(0, \infty)$ of compact support $[1, 2]$ with  $\phi^{(j)}(x)\ll 1$ for $j\geq 1$.
Let $X>1$ and $\a,\ \b\geq 0$.

{\rm (i)} Suppose
\begin{eqnarray}
2\max\{1,2^{\b-1/m}\}(\a\beta)^m\leq X^{1-\beta m}.
\end{eqnarray}
Then the estimate
\begin{eqnarray}
\sum_{n\neq0}A_f(n,1,\cdots,1)e(\pm\alpha |n|^{\beta})\phi\Big(\frac {|n|}X\Big)
\ll X^{-M}
\end{eqnarray}
holds for any $M>0$, where the implied constant in (1.12) depends  on $m$, $f$, $\beta$ and $M$.

{\rm (ii)} Suppose
\bea
2\max\{1,2^{\b-1/m}\}(\a\beta)^m> X^{1-\beta m}.
\eea
Then for $ \b\not=1/m$, we have
\begin{eqnarray*}
\sum_{n\neq0}A_f(n,1,\cdots,1)e(\pm\alpha |n|^{\beta})\phi\Big(\frac {|n|}X\Big)
\ll_{m,f,\b}(\a X^\b)^{m/2}.
\end{eqnarray*}
For $\b=1/m$, we have
\begin{eqnarray}
\sum_{n\neq0}A_f(n,1,\cdots,1)e(\pm\alpha |n|^{\beta})\phi\Big(\frac {|n|}X\Big)
\ll_{m,f}(\a X^\b)^{(m+1)/2}.
\end{eqnarray}
Moreover, for $X>\alpha^{m (m-1)/(1-m\varepsilon)}$ with $0<\varepsilon<1/m$,
\begin{eqnarray}
&&
\sum_{n\neq0}
A_f(n,1,\cdots,1)e\big(\pm \a |n|^{1/m} \big)
\phi\Big(\frac{|n|}{X}\Big)
\nonumber
\\
&=&m\le\{A_f(1,1,\cdots, 1,n_\a)+A_f(1,1,\cdots, 1,-n_\a)\ri\}
\sum_{k=0}^{r}\rho_{\pm}(k,m,\a, X )X^{ 1/(2m)+1/2-k/m}\nonumber
\\
&&+O_{m,r, \ve}\big(X^{-r/m+1/2+\varepsilon}\big),
\end{eqnarray}
where
\begin{eqnarray*}
\rho_{\pm}(k,m,\a, X)&=&(\mp i)^{k+(m-1)/2}\tilde{c}_kn_\a^{1/(2m)-1/2-k/m}\\
&&\times \int_0^\infty t^{m/2-k-1/2}\phi(t^m)e\big((\alpha- mn_\a^{1/m})X^{1/m}t\big)dt
\end{eqnarray*}
with $n_\a$ being the integer satisfying $(\a/m)^m-n_\a\in (-1/2,1/2]$,  $\tilde{c}_k$ being constants depending on $m$ and $f$.}

We remark that the condition (1.11) holds for any fixed $\a$ and $\b\in (0,1/m)$ when $X$ is large enough in terms of $\a$ and $\b$, and hence the rapid decay in (1.12) holds.
For $\b=1/m$ this is the case when $0\leq \a\leq m2^{-1/m} $. For $\b>1/m$, (1.12) holds for $0\leq \a\leq 2mX^{1/m-\b}$.
When $\a=0$ this result was firstly obtained by Booker \cite{Boo} via a different approach.

Note that the main term in (1.15) is negligible when $|(\alpha/m)^m-n_\a|>X^{\ve-1/m}$ since the integral in $\rho_{\pm}(k,m,\a, X)$ is arbitrarily small by repeated partial integrating by parts.
Consequently, (1.15) manifests resonance of Fourier coefficients of $f$ against
$e\big(\pm \alpha n^{1/m})$ when $\a$ approaches $mq^{1/m}$. In particular, when $\alpha=mq^{1/m}$ we obtain the following.
\medskip

\noindent
{\bf Corollary 1.1.} {\it Let $q$ be a positive integer and $0<\varepsilon<1/m$. Then for $X>(m^mq)^{(m-1)/(1-m\ve)}$ we have
\begin{eqnarray*}
&&
\sum_{n\neq0}
A_f(n,1,\ldots,1)e\big(\pm m(q |n|)^{1/{m}} \big)
\phi\Big(\frac{|n|}{X}\Big)
\nonumber
\\
&=&\le\{A_f(1,\ldots, 1,q) +A_f(1,\ldots, 1,-q) \ri\}
\sum_{k=0}^{r}\omega_{\pm}(k,m,q)X^{1/(2m)+1/2-k/m}\nonumber
\\
&&+O_{m, r, \ve}\big(X^{-r/m+1/2+\varepsilon}\big),
\end{eqnarray*}
where
\begin{eqnarray*}
\omega_{\pm}(k,m,q)
&=&(\mp i)^{k+(m-1)/2}\tilde{c}_kq^{1/(2m)-1/2-k/m}
\int_0^\infty x^{1/(2m)-1/2-k/m}\phi(x)dx.
\end{eqnarray*}}

The asymptotic behavior as described in Corollary 1.1 was proved for cusp forms for $SL_2(\mathbb Z)$ in
Iwaniec-Luo-Sarnak \cite{IwaLuoSar} and Ren-Ye \cite{Ren-Ye}, while for Maass forms
for $SL_3(\mathbb Z)$ it was proved in Ren-Ye \cite{RenYe2}. For similar
results concerning coefficients of $L$-functions for general Selberg class, see Kaczorowski and Perelli \cite{Per1}
\cite{Per2}.

We will also prove a sharp-cut version of Theorem 1.2.

\noindent
{\bf Theorem 1.3.} {\it Let $f$ be a full-level cusp form for
$GL_m(\mathbb Z)$. Let  $X>1$ and $\a, \ \b\geq 0$.

{\rm (i)} Assume that (1.11) holds. Then we have
\begin{eqnarray}
\sum_{X<|n|\leq 2X}A_f(n,1,\cdots,1)e(\pm\alpha |n|^{\beta})
\ll_{f,\b,m} X^{1-1/m}.
\end{eqnarray}

{\rm (ii)} Assume that (1.13) holds. Then for $\b\not=1/m$ we have
\begin{eqnarray}
\sum_{X<|n|\leq 2X}A_f(n,1,\cdots,1)e(\pm\alpha |n|^{\beta})
\ll_{f,\b,m} (\a X^\b)^{m/2}+X^{1-1/m}.
\end{eqnarray}
For $\b=1/m$ we have
\begin{eqnarray}
&&\sum_{X<|n|\leq 2X}A_f(n,1,\cdots,1)e(\pm\a |n|^{1/m})\nonumber\\
&&\ll \a^{1/2}X^{ 1/(2m)+1/2}+\a^{m-1/2}X^{1/2-1/(2m)}+X^{1-1/m}.
\end{eqnarray}}

\noindent
{\bf Theorem 1.4.} {\it Assume the following bound toward the Ramanujan conjecture
$$
A_f(1,\cdots, 1, n)\ll n^{\th},\qquad \rm{for\ some}\quad 0<\th<\frac2{m-1}.
$$

{\rm (i)} Suppose that the parameters $\a, \b, X$ satisfy (1.11). Then we have
\begin{eqnarray}
\sum_{X<|n|\leq 2X}A_f(n,1,\cdots,1)e(\a |n|^\b)
\ll_m X^{(m-1)(1+\th)/(m+1)}.
\end{eqnarray}

{\rm (ii)} Suppose that the parameters $\a, \b, X$ satisfy (1.13). Then we have, for $\b\not=1/m$,
\begin{eqnarray}
\sum_{X<|n|\leq 2X}A_f(n,1,\cdots,1)e(\a |n|^\b)\ll  (\a X^\b)^{m/2}+X^{(m-1)(1+\th)/(m+1)},
\end{eqnarray}
and for $\b=1/m$,
\begin{eqnarray}
&&\sum_{X<|n|\leq 2X}A_f(n,1,\cdots,1)e(\pm\a |n|^{1/m})\nonumber
\\
&&=-\sqrt{m}X^{ 1/(2m)+1/2}(-i)^{(m-1)/2}I(m,\a,X)\frac{A_f(1,\cdots, 1, n_\a)+A_f(1,\cdots, 1, -n_\a)}{
n_\a^{1/2-1/(2m)}}
\nonumber
\\
&&+O(\a^{m-1/2}X^{1/2-1/(2m)})+O(X^{(m-1)(1+\th)/(m+1)}),
\end{eqnarray}
where $n_\a$ is the integer satisfying $(\a/m)^m-n_\a\in(-1/2, 1/2]$ and
$$
I(m,\a,X)=\int_{1}^{2^{1/m}}t^{m/2-1/2}e\big((\alpha-mn_\a^{1/m} )X^{1/m}t\big)dt.
$$}

Note that for fixed $\alpha$ and large $X$, (1.21) is an asymptotic formula for $\theta<1/3$ when $m=3$ and for
$\theta<1/24$ when $m=4$. Recall
that the best known results are $\theta=5/14$ for $m=3$ and
$\theta=9/22$ for $m=4$, by Kim and Sarnak \cite{Kim}\cite{KimSar} and Sarnak \cite{Sar05} (21). When $m\geq 5$, (1.21) implies the following bound
\begin{eqnarray*}
&&\sum_{X<|n|\leq 2X}A_f(n,1,\cdots,1)e(\a |n|^{1/m})\nonumber
\\
&&\ll \a^{1/2+(\th-1/2)m}X^{ 1/(2m)+1/2}+\a^{m-1/2}X^{1/2-1/(2m)}+X^{(m-1)(1+\th)/(m+1)}.
\end{eqnarray*}
We point out that when $\th<1/m$, the bounds in (1.19)-(1.21) improves the bounds in (1.16)-(1.18), respectively.

We will prove Theorem 1.1 in Sections 2-4, and prove Theorems 1.2-1.4 in Sections 5-6.
\bigskip

\section{Stirling's asymptotic}
\setcounter{equation}{0}
\medskip

Changing $s$ to $1-2s$ and noting that
$\{{\mu}_{\tilde f}(j)\}=\{\overline\mu_f(j)\}$ for cusp form $f$,
we get
\begin{eqnarray}
\Psi(x)
=i\pi^{-{m}/2-1}
\int_{\text{Re}\, s=\s}
\big(\pi^mx\big)^{-2s+1}
\Tilde{\psi}(-2s+1)G(s)ds,
\end{eqnarray}
where
\begin{eqnarray}
G(s)=\prod_{j=1}^m
\frac{\Gamma\big(s-\frac{\overline\mu_f(j)}2\big)}
{\Gamma\big(-s+\frac{1-\mu_f(j)}2\big)}.
\end{eqnarray}
Note that the $\Gamma$-functions in the numerator on
the right side of (2.2) are analytic and nonzero for
\begin{eqnarray*}
\sigma>\frac12
\max\big\{{\rm Re}\,\overline\mu_f(1) ,\ldots,{\rm Re}\,\overline\mu_f(m)
\big\}.
\end{eqnarray*}
A bound due to Luo, Rudnick and Sarnak \cite{LRS} asserts that
\begin{eqnarray}
|{\rm Re}\,\mu_f(j)|\leq \frac12-\frac1{m^2+1},\quad j=1,\ldots,m.
\end{eqnarray}
Thus we are allowed to take any
$\sigma>1/4-1/(2(m^2+1))$ in (2.1) where the convergence of the integral is guaranteed by the
rapidly decay of $\Tilde{\psi}(-2s+1)$ with respect to $t.$ Let us consider $s=\s+it$ with
\begin{eqnarray}
\sigma>\frac14-\frac1{2(m^2+1)}.
\end{eqnarray}
Let $r>m/2$ and define
\begin{eqnarray}
\s(r)=\frac14+\frac r{2m}-\varepsilon,
\end{eqnarray}
where $\varepsilon>0$ is a small number.
We will show that in the vertical strip
$$
L_r:\ |{\rm Re}\, s-\s(r)|\leq \frac{\varepsilon}2,
$$
there holds
\begin{eqnarray}
G(s)=m^{-2ms+m/2}\cdot
\frac{\Gamma\big(ms-\frac{m-1}2\big)}
{\Gamma\big(-ms+\frac12\big)}\cdot
\Big(1+\sum_{j=1}^r\frac{h_j}{s^{j}}+ E_r(s) \Big),
\end{eqnarray}
where $h_j$'s are constants depending on $f$, $E_r(s)$ is analytic in $L_r$ and satisfies $E_r(s)=O(|s|^{-r-1})$ as $|s|\to \infty$.

By Stirling's formula \cite{RenYe4} \cite{Spi}, for $|{\rm Im} s|\geq 2|\b|$ and $|{\rm Im} s|\gg |s|$ one has
\begin{eqnarray}
\log \Gamma(s+\b)
=\Big(s+\b-\frac12\Big)\log s-s
+\log\sqrt{2\pi}+\sum_{j=1}^r\frac{d_j}{s^{j}}
+O\big(|s|^{-r-1}\big),
\end{eqnarray}
where $d_j$ are constants depending on $\b,$ and the implied constant depends on $\b$ and $r$.
This in combination with the fact that $\mu_f(1)+\cdots+\mu_f(m)=0$ show that, for $s$ satisfying
\begin{eqnarray}
|\mbox{Im}\ s|\gg |s|\quad {\rm and}\quad |\mbox{Im}\ s|\geq t_0=2+\max_{1\leq j\leq m}\big\{|\mu_f(j)|\big\},
\end{eqnarray}
there holds
\begin{eqnarray}
\log G(s)
&=&
\Big(ms-\frac{m}2\Big)\log s-2ms+ms\log (-s)
\nonumber
\\
&+&
\sum_{j=1}^r\frac{f_j}{s^{j}}+O\big(|s|^{-r-1}\big),
\end{eqnarray}
where $f_j$ are constants depending on $m$ and $\mu_f(j)$, $j=1,\cdots, m$.
Similarly, by (2.7), for $|\mbox{Im}\, s|\geq 1$, we have
\begin{eqnarray}
\log \frac{\Gamma\big(ms-\frac{m-1}2\big)}
{\Gamma\big(-ms+\frac12\big)}
&=&
\Big(ms-\frac{m}2\Big)\log ms-2ms+ms\log (-ms)
\nonumber
\\
&+&
\sum_{j=1}^r\frac{g_j}{s^{j}}+O\big(|s|^{-r-1}\big).
\end{eqnarray}
Comparing (2.9) with (2.10) we get, for $s$ satisfying (2.8), that
\begin{eqnarray}
G(s)=m^{-2ms+m/2}\cdot
\frac{\Gamma\big(ms-\frac{m-1}2\big)}
{\Gamma\big(-ms+\frac12\big)}\cdot
\Big(1+\sum_{j=1}^r\frac{h_j}{s^{j}}
+ O\big(|s|^{-r-1}\big) \Big).
\end{eqnarray}

On the other hand, one can write
\begin{eqnarray}
G(s)
=m^{-2ms+m/2}\cdot
\frac{\Gamma\big(ms-\frac{m-1}2\big)}
{\Gamma\big(-ms+\frac12\big)}
\cdot A(s),
\end{eqnarray}
where $A(s)$ is analytic and non-zero for $s=\s+it$ satisfying (2.4) and
\begin{eqnarray}
&&
-ms+\frac{m-1}2\not\in\mathbb N,
\quad
ms-\frac12\not\in\mathbb N,
\nonumber
\\
&&
-\Big(s-\frac{\overline\mu_f(j)}2\Big)\not\in\mathbb N,
\quad
s-\frac{1-\mu_f(j)}2\not\in\mathbb N,
\quad
1\leq j\leq m.
\end{eqnarray}
Here we have to avoid poles and zeros of the quotients of $\Gamma$-functions on both sides of (2.12). Let $\s(r)$ be defined by (2.5).
We show that when $\varepsilon>0$ is small enough (2.13) is satisfied for
$s=\s+it$ with $|\s-\s(r)|\leq \varepsilon/2$. Actually, let $\|a\|$ denote the distance from $a$ to the nearest integer.
One can easily check that  $-ms+(m-1)/2\not\in\mathbb N$,
$ ms-1/2\not\in\mathbb N$ and
$-(s-\overline\mu_f(j)/2)\not\in\mathbb N$ when $\varepsilon<1/(8m)$ and $|\s-\s(r)|\leq\varepsilon/2$.
Moreover, one has
\begin{eqnarray}
\Big\|\s(r)-\frac{1-{\rm Re} \mu_f(j)}2\Big\|
&=&
\Big\|\frac{r}{2m}-\frac14+\frac{{\rm Re} \mu_f(j)}2-\varepsilon\Big\|.
\end{eqnarray}
If $m|r$, the right expression in (2.14) is $\geq 1/{2(m^2+1)}-\varepsilon>\varepsilon/2$ when $\varepsilon<1/{4(m^2+1)}$. For $m\nmid r$, write
$$
\d_j=\Big\|\frac{r}{2m}-\frac14+\frac{{\rm Re} \mu_f(j)}2\Big\|.
$$
Note that when $\d_j=0$, the right expression in (2.14) is $\geq \varepsilon$ for any $0<\varepsilon<1/2$. If $\d_j\not=0$, the expression is
$>\varepsilon$ when $\varepsilon<\d_j/2.$ Thus one can choose
$$
0<\varepsilon<\min_{
\genfrac{}{}{0pt}{}{1\leq j\leq m}{\d_j\not=0}
}
\Big\{\frac1{4(m^2+1)},\ \frac{\d_j}2\Big\}
$$
so that (2.13) is satisfied for $|\s-\s(r)|\leq \varepsilon/2.$

\medskip

Now we consider the region
$$
D=\Big\{s=\s+it\ \Big|\ |\s-\s(r)|\leq \frac{\varepsilon}2,\ |t|\leq |\s(r)|+2t_0\Big\}.
$$
In this region one can express $A(s)$ as
\begin{eqnarray*}
1+\sum_{j=1}^r\frac{h_j}{s^{j}}+ E_r(s),
\end{eqnarray*}
where $E_r(s)$ is analytic in $D$.
Back to (2.12) and (2.11), we see that (2.6) holds in the vertical strip
$L_r=\{\s+it: |\s-\s(r)|\leq {\varepsilon}/2\}$.

\medskip

Write $u=\pi^mx$ and define
\begin{eqnarray}
\mathcal {H}_{j}
&=&
i\pi^{-m/2-1}
\int_{\text{Re}\, s=\s(r)}
\frac{\Gamma\big(ms-\frac{m-1}2\big)}
{s^{j}\Gamma\big(-ms+\frac12\big)}
\nonumber
\\
&\times&
m^{-2ms+m/2}u^{-2s+1}
\Tilde{\psi}(-2s+1)ds,
\\
\mathcal {E}_r
&=&
i\pi^{-m/2-1}
\int_{\text{Re}\, s=\s(r)}
\frac{m^{-2ms+m/2}\Gamma\big(ms-\frac{m-1}2\big)}
{\Gamma\big(-ms+\frac12\big)}
\nonumber
\\
&\times&
E_r(s)u^{-2s+1}\Tilde{\psi}(-2s+1)ds.
\end{eqnarray}
Then we get
\begin{eqnarray*}
\Psi(x)=
\mathcal {H}_{0}+
\sum_{j=1}^rh_j\mathcal {H}_{j}+\mathcal {E}_r.
\end{eqnarray*}
We remark that the integral for $\mathcal {E}_r$ is only valid in
the vertical strip $L_r$, while the integral path for $\mathcal {H}_{j}$
can be moved freely in the half plane
$\text{Re}\, s> 1/2-1/(2m)$. In the following we will replace $\s(r)$ by $\s_1(r)$ in (2.15) where
\begin{eqnarray}
\s_1(r)=\frac12-\frac1{2m}+\frac{r}{m}+\varepsilon.
\end{eqnarray}

The estimate of $\mathcal {E}_r$ is immediate. By Stirling's
formula, for  $s=\s+it$ we have
\begin{eqnarray*}
\frac{\Gamma\big(ms-\frac{m-1}2\big)}
{\Gamma\big(-ms+\frac12\big)}
\ll |s|^{2m\s-m/2},\quad {\rm as }\ |s|\to\infty.
\end{eqnarray*}
When $r>m/2$, $\s=\s(r)$ as defined in (2.5) one has
\begin{eqnarray*}
\frac{\Gamma\big(ms-\frac{m-1}2\big)}
{\Gamma\big(-ms+\frac12\big)}E_r(s)
\ll |s|^{2m\s(r)-m/2-r-1}\ll |s|^{-1-\varepsilon}.
\end{eqnarray*}
Note that
\begin{eqnarray}
\tilde\psi(s)= \int_{aX}^{bX} \psi(y)y^{s}\frac {dy}y\ll X^{\s}\int_a^b|\phi(x)|x^{\s-1}dx\ll_{\phi} X^\s.
\end{eqnarray}
Hence
\begin{eqnarray*}
\mathcal {E}_r&\ll& u^{-2\s(r)+1}X^{-2\s(r)+1}
\int_{\text{Re}\, s=\s(r)}|s|^{-1-\varepsilon}|ds|
\nonumber
\\
&\ll&(uX)^{-2\s(r)+1}\ll(xX)^{-{r}/m+1/2+2\varepsilon}.
\end{eqnarray*}
This shows
\begin{eqnarray}
\Psi(x)=
\mathcal {H}_{0}+
\sum_{j=1}^rh_j\mathcal {H}_{j}
+O\big((xX)^{-{r}/m+1/2+\varepsilon}\big).
\end{eqnarray}

To prove Theorem 1.1 it remains to compute $\mathcal {H}_{j}$ for $j=0, 1, \ldots , r$, this will be carried out in the following two
sections according to $m$ is odd or even.

\bigskip

\section{Proof of Theorem 1.1 when $m$ is odd}
\setcounter{equation}{0}

\medskip

By definition we have
\begin{eqnarray*}
\mathcal {H}_{0}
=i\pi^{-m/2-1}\int_{\text{Re}\, s=\s_1(r)}
\frac{\Gamma\big(ms-\frac{m-1}2\big)}
{\Gamma\big(-ms+\frac12\big)}
m^{-2ms+m/2}u^{-2s+1}
\Tilde{\psi}(-2s+1)ds.
\end{eqnarray*}
Move the integral path to ${\rm Re}\, s=-\infty$ and
note that the only poles of the integrand are from $\Gamma(ms-(m-1)/2)$ and hence are simple at
$s=-n/m+(m-1)/(2m)$  for $n=0, 1, \ldots $ with the residue $(-1)^n/(n!m)$. Thus
\begin{eqnarray*}
\mathcal {H}_{0}
&=&
-2\pi^{-m/2}\sum_{n=0}^\infty
\frac{(-1)^n}{n!\Gamma\big(n+1-\frac{m}2\big)}
m^{2n-m/2}u^{(2n+1)/m}
\Tilde{\psi}\Big(\frac{2n+1}m\Big)
\\
&=&
-2\pi^{-m/2}
\int_0^\infty
\sum_{n=0}^\infty
\frac{(-1)^n}{n!\Gamma\big(n+1-\frac{m}2\big)}
m^{2n-m/2}
(uy)^{(2n+1)/m}
\psi(y)
\frac{dy}y.
\end{eqnarray*}
Recall the power series definition of the Bessel function
\begin{equation}
J_\nu(z)
=
\sum_{n=0}^\infty
\frac{(-1)^n}{n!\Gamma\big(n+1+\nu\big)}
\Big(\frac z2\Big)^{2n+\nu},
\end{equation}
we get
\begin{eqnarray}
\mathcal {H}_{0}
&=&
-2\pi^{-m/2}u
\int_0^\infty
\sum_{n=0}^\infty
\frac{(-1)^n}{n!\Gamma\big(n+1-\frac{m}2\big)}
\nonumber
\\
&\times&
(uy)^{1/m-1/2}
\big(m(uy)^{1/m}\big)^{2n-m/2}
\psi(y)dy
\nonumber
\\
&=&
-2\pi^{-m/2}u
\int_0^\infty
(uy)^{1/m-1/2}
J_{-m/2}\big(2m(uy)^{1/m}\big)
\psi(y)
dy.
\end{eqnarray}

Note that if $\nu=-q-1/2$ and $q\geq 1$, the half
integral-order Bessel function is elementary. Applying
(8.462.2) in \cite{GraRyz}, we have
\begin{eqnarray}
J_{-q-1/2}(z)
&=&
\frac1{\sqrt{2\pi z}}\sum_{j=0}^{q}
\frac{(q+j)!}{j!(q-j)!(2z)^{j}}
\nonumber
\\
&\times&
\Big\{
i^{j+q}e\Big(\frac{z}{2\pi}\Big)
+(-i)^{j+q}e\Big(-\frac{z}{2\pi}\Big)
\Big\}.
\end{eqnarray}
When $m=2\ell+1$, we apply (3.3) with
$z=2m(uy)^{1/m}=2\pi m(xy)^{1/m} $, $q=\ell$ and then substitute it back to (3.2) to get
\begin{eqnarray*}
\mathcal H_0
&=&
-2\pi^{m/2}x
\int_0^\infty
\big(\pi^mxy\big)^{1/m-1/2}
\frac1{\sqrt{4\pi^2m(xy)^{1/m}}}
\\
&\times&
\sum_{j=0}^\ell
\frac{(\ell+j)!}{j!(\ell-j)!\big(4\pi m(xy)^{1/m}\big)^{j}}
\\
&\times&
\Big\{
i^{j+q}e\big(m(xy)^{1/m} \big)
+(-i)^{j+q}e\big(-m(xy)^{1/m} \big)
\Big\}
\psi(y)dy.
\end{eqnarray*}
Denoting
\begin{eqnarray*}
a_j^{(0)}
=
-\frac{1}{\sqrt{m}}\cdot
\frac{(\ell+j)!}{j!(\ell-j)!(4\pi m)^{j}},
\quad j=0,\ 1,\ \ldots ,\ \ell,
\end{eqnarray*}
we get
\begin{eqnarray}
\mathcal {H}_{0}
&=&
x\sum_{j=0}^{\ell}
a_j^{(0)}
\int_0^\infty
(xy)^{-(j+\ell)/m}\psi(y)
\nonumber
\\
&\times&
\Big\{i^{j+\ell}e\big(m(xy)^{1/m}\big)
+(-i)^{j+\ell}e\big(-m(xy)^{1/m}\big)\Big\}dy.
\end{eqnarray}

Now we turn to the estimate of $\mathcal {H}_{1}$. We have
\begin{eqnarray*}
\mathcal {H}_{1}
&=&
i\pi^{-m/2-1}
\int_{\text{Re}\, s=\s_1(r)}
\frac{\Gamma(ms-\ell)}{s\Gamma\big(-ms+\frac12\big)}
\\
&\times&
m^{-2ms+m/2}u^{-2s+1}
\Tilde{\psi}(-2s+1)ds.
\end{eqnarray*}
Using the formula $\Gamma(s+1)=s\Gamma(s)$ repeatedly we obtain
\begin{eqnarray*}
\frac{\Gamma(ms-\ell)}{s}
&=&
m\Big\{\sum_{b=0}^{r-1}\frac{(\ell+b)!}{\ell!}(-1)^{b}\Gamma(ms-(\ell+b+1))\Big\}\\
&+&
(-1)^r(\ell+1)\cdots (\ell+r)\frac{\Gamma(ms-(\ell+r))}{s}.
\end{eqnarray*}
This gives us
\begin{eqnarray}
\mathcal {H}_{1}
&=&
i\pi^{-m/2-1}\Big(\sum_{b=0}^{r-1}\frac{(\ell+b)!}{\ell!}(-1)^{b}m\mathcal{I}_{b+1}\Big)
\nonumber
\\
&+&
i(-1)^r\pi^{-m/2-1}(\ell+1)\cdots (\ell+r)\mathcal{I}_{r}^*,
\end{eqnarray}
where for $d=1, ..., r$
\begin{equation}
\mathcal{I}_d=\int_{\text{Re}\, s=\s_1(r)}
\frac{\Gamma(ms-(\ell+d))}{\Gamma\big(-ms+\frac12\big)}
m^{-2ms+m/2}u^{-2s+1}
\Tilde{\psi}(-2s+1)ds,
\end{equation}
and
\begin{equation}
\mathcal{I}_r^*=\int_{\text{Re}\, s=\s_1(r)}
\frac{\Gamma(ms-(\ell+r))}{s\Gamma\big(-ms+\frac12\big)}
m^{-2ms+m/2}u^{-2s+1}
\Tilde{\psi}(-2s+1)ds.
\end{equation}

Moving the integral contour in (3.6) to Re$(s)=-\infty$ and picking up residues $(-1)^n/(mn!)$ at
$s=(-n+\ell+d)/m$ for $n=0,\ 1,\ 2, ...$, we get
\begin{eqnarray*}
m\mathcal{I}_d
&=&
2\pi i
\sum_{n=0}^{\infty}
\frac{(-1)^n}{n!\Gamma\big(n-(\ell+d)+\frac12\big)}
\\
&\times&
m^{2(n-\ell-d)+m/2}u^{2(n-\ell-d)/m+1}
\Tilde{\psi}\Big(\frac{2(n-\ell-d)}m+1\Big)
\\
&=&
2\pi i\sum_{n=0}^{\infty}\frac{(-1)^n}{n!\Gamma\big(n-(\ell+d)+\frac12\big)}
\\
&\times&
m^{2(n-\ell-d)+m/2}u^{2(n-\ell-d)/m+1}
\int_0^\infty
y^{2(n-\ell-d)/m}
\psi(y)dy.
\end{eqnarray*}
Collecting factors and using (3.1), we have
\begin{eqnarray*}
m\mathcal{I}_d
&=&
2\pi i um^{m/2+1/2-(l+d)}
\int_0^\infty
(uy)^{1/(2m)-(\ell+d)/m}
\\
&\times&
\sum_{n=0}^\infty
\frac{(-1)^n}{n!\Gamma\big(n-(\ell+d)+\frac12\big)}
\big(
m(uy)^{1/m}
\big)^{2n-\ell-d-1/2}
\psi(y)dy
\\
&=&
2\pi i um^{m/2+1/2-(l+d)}
\\
&\times&
\int_0^\infty
(uy)^{1/(2m)-(\ell+d)/m}
J_{-\ell-d-1/2}\big(2m(uy)^{1/m}\big)
\psi(y)dy.
\end{eqnarray*}
Using (3.3) with $z=2m(uy)^{1/m}=2\pi m(xy)^{1/m}$ and $q=\ell+d$ and then following the computation leading to
(3.4), we get
\begin{eqnarray}
m\mathcal{I}_d
&=&
i x\pi^{m/2+1-d}
m^{1/2-d}
\sum_{j=0}^{\ell+d}
\frac{(\ell+d+j)!}{j!(\ell+d-j)!(4\pi m)^j}
\int_0^\infty(xy)^{-(\ell+d+j)/m}
\nonumber
\\
&\times&
\Big\{i^{\ell+d+j}e\big(m(xy)^{1/m}\big)
+(-i)^{\ell+d+j}e\big(-m(xy)^{1/m}\big)\Big\}
\psi(y)dy.
\end{eqnarray}

To estimate $\mathcal{I}_r^*,$ we move the integral contour
in (3.7) from Re$(s)=\sigma_1(r)$ back to
Re$(s)=\sigma(r)=1/4+r/(2m)-\varepsilon$. Then the real part of
$ms-(\ell+r)$ moves from $\varepsilon$ to
$1/4-(\ell+r)/2-m\varepsilon$. We pick up residues at
$s=(-n+\ell+r)/m$ for
$n=0,\ldots,[(\ell+r)/2-1/4]$ to get
\begin{eqnarray}
\mathcal{I}_r^*
&=&2\pi i um^{m/2}
\sum_{n=0}^{[(\ell+r)/2-1/4]}
\frac{(-1)^n}{(\ell+r-n)n!\Gamma\big(n-\ell-r+\frac12\big)}
\nonumber
\\
&\times&
\int_0^\infty
\big(m(uy)^{1/m}\big)^{2n-2(\ell+r)}
\psi(y)dy
\nonumber
\\
&+&
\int_{\text{Re}\, s=\s(r)}
\frac{\Gamma(ms-(\ell+r))}{s\Gamma\big(-ms+\frac12\big)}
m^{-2ms+m/2}u^{-2s+1}
\Tilde{\psi}(-2s+1)ds.
\end{eqnarray}
Since $\psi(y)(=\phi(y/X))$ is supported on $[aX,bX]$, the first term on the right side above is
\begin{eqnarray*}
\ll (uX)^{1+\frac{2[(\ell+r)/2-1/4]-2(\ell+r)}{m}}
\ll (xX)^{1-{(\ell+r)}/m-1/(2m)}
=(xX)^{-{r}/m+1/2},
\end{eqnarray*}
where the implied constant depends on $\phi$, $a$ and $b$. By Stirling's formula, for Re$(s)=\s(r)$ one has
\begin{eqnarray*}
\frac{\Gamma(ms-(\ell+r))}{s\Gamma\big(-ms+\frac12\big)}
\ll
|s|^{2m\s(r)-\ell-r-\frac32}
\ll
|s|^{-1-\varepsilon}.
\end{eqnarray*}
By (2.18) the last integral in (3.9) is
\begin{eqnarray*}
\ll (uX)^{-2\s(r)+1}\ll (xX)^{-{r}/m+1/2+\varepsilon}.
\end{eqnarray*}
Therefore
\begin{equation}
\mathcal I^*_r
\ll
(uX)^{-2\s(r)+1}\ll (xX)^{-{r}/m+1/2+\varepsilon}.
\end{equation}
This together with (3.8) and (3.5) show that
\begin{eqnarray*}
\mathcal {H}_{1}
&=&
-xm^{1/2}
\sum_{b=0}^{r-1}
(\pi m)^{-b-1}\frac{(\ell+b)!}{\ell!}
(-1)^{b}
\\
&\times&
\sum_{j=0}^{\ell+b+1}
\frac{(\ell+b+1+j)!}{j!(\ell+b+1-j)!(4\pi m)^j}
\int_0^\infty
(xy)^{-(\ell+b+1+j)/m}
\\
&\times&
\Big\{i^{\ell+b+1+j}e\big(m(xy)^{1/m}\big)
+(-i)^{\ell+b+1+j}e\big(-m(xy)^{1/m}\big)\Big\}
\psi(y)dy
\\
&+&
O\big((xX)^{-{r}/m+1/2+\varepsilon}\big).
\end{eqnarray*}
Collecting like terms, we get
\begin{eqnarray*}
\mathcal {H}_{1}
&=&
x\sum_{t=1}^{\ell+2r}
a_t^{(1)}
\int_0^\infty
(xy)^{-(\ell+t)/m}
\\
&\times&
\Big\{i^{\ell+t}e\big(m(xy)^{1/m}\big)
+(-i)^{\ell+t}e\big(-m(xy)^{1/m}\big)\Big\}
\psi(y)dy
\\
&+&
O\big((xX)^{-{r}/m+1/2+\varepsilon}\big),
\end{eqnarray*}
where for $t=1,\ 2, \ldots, \ell+2r,$
\begin{eqnarray*}
a_t^{(1)}
=-\frac{4\sqrt m}{(4\pi m)^t}\cdot
\frac{(\ell+t)!}{\ell!}
\sum_{\max\big\{0,\frac{t-\ell}{2}-1\big\}
\leq b\leq \min(r-1,t-1)}
\frac{(-4)^{b}(\ell+b)!}{(t-b-1)!(\ell+2b-t+2)!}.
\end{eqnarray*}
Note that when $j\geq r+1,$
\begin{eqnarray*}
&&
x\int_0^\infty(xy)^{-(\ell+j)/m}
\Big\{i^{\ell+j}e\big(m(xy)^{1/m}\big)
+(-i)^{\ell+j}e\big(-m(xy)^{1/m}\big)\Big\}
\psi(y)dy
\\
&\ll&
(xX)^{1-{(\ell+r+1)}/m}=(xX)^{-{r}/m+1/2-1/(2m)}.
\end{eqnarray*}
Thus we finally obtain
\begin{eqnarray}
\mathcal {H}_{1}
&=&
x\sum_{t=1}^{r}
a_t^{(1)}
\int_0^\infty(xy)^{-(\ell+t)/m}
\nonumber
\\
&\times&
\Big\{i^{\ell+t}e\big(m(xy)^{1/m}\big)
+(-i)^{\ell+t}e\big(-m(xy)^{1/m}\big)\Big\}
\psi(y)dy
\nonumber
\\
&+&
O\big((xX)^{-{r}/m+1/2+\varepsilon}\big).
\end{eqnarray}

To finish the proof of Theorem 1.1, we need to estimate
$\mathcal {H}_{j}$ for $2\leq j\leq r.$
This is similar to the case when $j=1$ and so we will briefly describe the idea. By (2.15) we have
\begin{eqnarray*}
\mathcal {H}_{j}
=i\pi^{-m/2-1}\int_{\text{Re}\, s=\s_1(r)}
\frac{\Gamma(ms-\ell)}{s^j\Gamma\big(-ms+\frac12\big)}
m^{-2ms+m/2}u^{-2s+1}
\Tilde{\psi}(-2s+1)ds.
\end{eqnarray*}
By repeated use of the formula $\Gamma(s+1)=s\Gamma(s),$ we have
\begin{eqnarray*}
\frac{\Gamma(ms-\ell)}{s^j}
&=&
m\Big\{
\frac{\Gamma(ms-(\ell+1))}{s^{j-1}}
\\
&+&
\sum_{k=1}^{h-1}(\ell+1)\cdots (\ell+k)(-1)^{k}
\frac{\Gamma(ms-(\ell+k+1))}{s^{j-1}}\Big\}
\\
&+&
(-1)^h(\ell+1)\cdots (\ell+h)\frac{\Gamma(ms-(\ell+h))}{s^j},
\end{eqnarray*}
where $h=r-j+1$. Applying this process repeatedly, we can
finally decompose
$\Gamma(ms-\ell)/(s^j\Gamma(-ms+1/2))$ into finite sums of
\begin{eqnarray*}
\lambda\frac{\Gamma(ms-(\ell+p))}{\Gamma\big(-ms+\frac12\big)},\quad
2\leq p\leq r,
\end{eqnarray*}
and
\begin{eqnarray*}
\mu\frac{\Gamma(ms-(\ell+p))}{s^q\Gamma\big(-ms+\frac12\big)},\quad
p+q=r+1,\ \quad 2\leq p\leq r, \quad 1\leq q\leq r
\end{eqnarray*}
with some constants $\lambda=\lambda_p, \mu=\mu_{p,q}$.
So it remains to estimate integrals of the form
\begin{eqnarray*}
\mathcal{I}_p=\int_{\text{Re}\, s=\s_1(r)}
\frac{\Gamma(ms-(\ell+p))}{\Gamma\big(-ms+\frac12\big)}
m^{-2ms+m/2}u^{-2s+1}
\Tilde{\psi}(-2s+1)ds,  \ \  2\leq p\leq r,
\end{eqnarray*}
and
\begin{eqnarray*}
\mathcal{I}^*=\int_{\text{Re}\, s=\s_1(r)}
\frac{\Gamma(ms-(\ell+p))}{s^q\Gamma\big(-ms+\frac12\big)}
m^{-2ms}u^{-2s}
\Tilde{\psi}(-2s+1)ds
\end{eqnarray*}
with $p+q=r+1,\ \  2\leq p\leq r, \ 1\leq q\leq r.$

\medskip

The integral $\mathcal{I}_p$ has been treated in (3.8).  The integral $\mathcal{I}^*$ can be estimated in a similar way as $\mathcal{I}_r^*$, and it
satisfies $\mathcal{I}^*\ll (uX)^{-r/m+1/2+\varepsilon}.$ So one finally  obtains for $j\geq 2$,
\begin{eqnarray}
\mathcal {H}_{j}
&=&
x\sum_{t=j}^{r}
a_t^{(j)}
\int_0^\infty(xy)^{-(\ell+t)/m}
\nonumber
\\
&\times&
\Big\{
i^{\ell+t}e\big(m(xy)^{1/m}\big)
+(-i)^{\ell+t}e\big(-m(xy)^{1/m}\big)
\Big\}
\psi(y)dy
\nonumber
\\
&+&
O\big((xX)^{-{r}/m+1/2+\varepsilon}\big).
\end{eqnarray}

Collecting estimates in (3.4), (3.11) and (3.12), and
substituting them into (2.19), we finishes the proof of Theorem 1.1 when $m$ is odd by writing $c_k=\sum_{j=0}^ka_k^{(j)}$. In particular we have $c_0=a_0^{(0)}=-1/\sqrt{m}$.

\bigskip

\section{Proof of Theorem 1.1 when $m$ is even}
\setcounter{equation}{0}

\medskip

By (3.2), when $m=2l$ one has
\begin{equation}
\mathcal {H}_{0}
=
-2\pi x
\int_0^\infty
(xy)^{1/m-1/2}
\psi(y)
J_{-\ell}\big(2m\pi (xy)^{1/m}\big)dy.
\end{equation}
Note that for positive integer $k,$ $J_{-k}(z)=(-1)^{k}J_{k}(z)$. By (4.8) in \cite{LiuYe1},
\begin{eqnarray*}
J_{k}(z)
&=&
\frac{1}{\sqrt{2\pi z}}e^{i(z-(2k+1)\pi/4)}
\sum_{0\leq j<2L}
\frac{i^j\Gamma\big(k+j+\frac12\big)}
{j!\Gamma\big(k-j+\frac12\big)(2z)^{j}}(2z)^{-j}
\nonumber
\\
&+&
\frac{1}{\sqrt{2\pi z}}e^{-i(z-(2k+1)\pi/4)} \sum_{0\leq j<2L}
\frac{i^j\Gamma\big(k+j+\frac12\big)}
{j!\Gamma\big(k-j+\frac12\big)(2z)^{j}}(-2z)^{-j}
\nonumber
\\
&+&
O\big(|z|^{-2L-1/2}\big).
\end{eqnarray*}
Since $e^{-i(2k+1)\pi/4}=(-i)^ki^{-1/2}$ and
$ e^{i(2k+1)\pi/4}=i^k(-i)^{-1/2}$ the above formula can be rewritten as
\begin{eqnarray}
J_{-k}(z)
&=&
\frac{1}{\sqrt{2\pi z}}
\sum_{0\leq j<2L}
\frac{\Gamma\big(k+j+\frac12\big)}
{j!\Gamma\big(k-j+\frac12\big)(2z)^{j}}
\nonumber
\\
&\times&
\Big\{i^{j+k-1/2}e\Big(\frac{z}{2\pi}\Big)
+(-i)^{j+k-1/2}e\Big(-\frac{z}{2\pi}\Big)\Big\}
\nonumber
\\
&+&
O\big(|z|^{-2L-1/2}\big).
\end{eqnarray}
Let $z=2m\pi (xy)^{1/m}$ and $k=\ell$ in (4.2), we get
\begin{eqnarray*}
\mathcal {H}_{0}
&=&
-\frac{x}{\sqrt{m}}
\sum_{0\leq j<2L}
\frac{\Gamma\big(\ell+j+\frac12\big)}
{j!\Gamma\big(\ell-j+\frac12\big)(4\pi m)^j}
\int_0^\infty
(xy)^{1/(2m)-(j+\ell)/m}
\psi(y)
\\
&\times&
\Big\{i^{j+\ell-1/2}e\big(m(xy)^{1/m}\big)
+(-i)^{j+\ell-1/2}e\big(-m(xy)^{1/m}\big)\Big\}
dy
\\
&+&
O\big((xX)^{1/2-(2L)/m+1/(2m)}\big).
\end{eqnarray*}
For any $r\geq 2$ and  $L=[r/2]+1$, the error term above is
$O((xX)^{1/2-r/m})$, and this gives
\begin{eqnarray*}
\mathcal {H}_{0}
&=&
x\sum_{j=0}^{2[r/2]+1}
b_j^{(0)}
\int_0^\infty
(xy)^{1/(2m)-(j+\ell)/m}
\psi(y)
\\
&\times&
\Big\{i^{j+\ell-1/2}e\big(m(xy)^{1/m}\big)
+(-i)^{j+\ell-1/2}e\big(-m(xy)^{1/m}\big)\Big\}dy
\\
&+&
O\big((xX)^{1/2-r/m}\big),
\end{eqnarray*}
where
\begin{eqnarray*}
b_j^{(0)}=-\frac{1}{\sqrt{m}}\cdot
\frac{\Gamma\big(\ell+j+\frac12\big)}
{j!\Gamma\big(\ell-j+\frac12\big)(4\pi m)^j}.
\end{eqnarray*}

Now we turn to the estimate of $\mathcal {H}_{j}$ for $j\geq 1$. By (2.15) we have
\begin{eqnarray*}
&&
\mathcal {H}_{j}
=i\pi^{-\ell-1}\int_{\text{Re}\, s=\sigma_1(r)}
\frac{\Gamma\big(ms-\big(\ell-\frac12\big)\big)}
{s^{j}\Gamma\big(-ms+\frac12\big)}
m^{-2ms+m/2}u^{-2s+1}
\Tilde{\psi}(-2s+1)ds.
\end{eqnarray*}
We first consider $\mathcal {H}_{1}$. Using the formula
$\Gamma(s+1)=s\Gamma(s),$ one obtains
\begin{eqnarray*}
\frac{\Gamma\big(ms-\big(\ell-\frac12\big)\big)}{s}
=
m\Gamma\Big(ms-\Big(\ell+\frac12\Big)\Big)
-\Big(\ell+\frac12\Big)
\frac{\Gamma\big(ms-\big(\ell+\frac12\big)\big)}{s}.
\end{eqnarray*}
Applying this formula to the quotient on its right
side repeatedly, we get
\begin{eqnarray*}
&&
\frac{\Gamma\big(ms-\big(\ell-\frac12\big)\big)}{s}
\\
&=&
m\Big\{
\Gamma\Big(ms-\Big(\ell+\frac12\Big)\Big)
\\
&+&
\sum_{b=1}^{r-1}(-1)^{b}
\Big(\ell+\frac12\Big)\cdots
\Big(\ell+\frac{2b-1}2\Big)
\Gamma\Big(ms-\Big(\ell+\frac{2b+1}2\Big)\Big)
\Big\}
\\
&+&
(-1)^r
\Big(\ell+\frac12\Big)\cdots
\Big(\ell+\frac{2r-1}2\Big)
\frac{\Gamma\big(ms-\big(\ell+\frac{2r-1}2\big)\big)}{s}.
\end{eqnarray*}
Substituting this into the integral defining
$\mathcal H_1$, we have
\begin{eqnarray}
\mathcal {H}_{1}
&=&
\frac{i\pi^{-\ell-1}}{\Gamma\big(\ell+\frac12\big)}
\sum_{b=0}^{r-1}(-1)^{b}
\Gamma\Big(\ell+b+\frac12\Big)
m\mathfrak{I}_{b+1}
\nonumber
\\
&+&
i\pi^{-\ell-1}(-1)^r
\Big(\ell+\frac12\Big)\cdots
\Big(\ell+r-\frac12\Big)
\mathfrak{I}_{r}^*,
\end{eqnarray}
where for $d=1,\ldots,r$
\begin{eqnarray*}
\mathfrak{I}_d=\int_{\text{Re}\, s=\s_1(r)}
\frac{\Gamma\big(ms-\big(\ell+d-\frac12\big)\big)}
{\Gamma\big(-ms+\frac12\big)}
m^{-2ms+m/2}u^{-2s+1}
\Tilde{\psi}(-2s+1)ds,
\end{eqnarray*}
and
\begin{eqnarray*}
\mathfrak{I}_r^*=\int_{\text{Re}\, s=\s_1(r)}
\frac{\Gamma\big(ms-\big(\ell+r-\frac12\big)\big)}
{s\Gamma\big(-ms+\frac12\big)}
m^{-2ms+m/2}u^{-2s+1}
\Tilde{\psi}(-2s+1)ds.
\end{eqnarray*}

For $\mathfrak{I}_d$ we move the integral line to
$\Re s=-\infty$ and pick up residues
$(-1)^n/(n!m)$ at $s=(-n+\ell+d-1/2)/m$
with $n=0,\ 1,\ 2,\ldots$. Consequently,
\begin{eqnarray*}
m\mathfrak{I}_d
&=&
2\pi i\sum_{n=0}^{\infty}
\frac{(-1)^nm^{2(n-\ell-d)+m/2+1}}
{n!\Gamma\big(n+1-\ell-d)\big)}
\\
&\times&
u^{2(n-\ell-d)/m+1/{m}+1}
\Tilde{\psi}\Big(\frac{2(n-\ell-d)}m+\frac1m+1\Big)
\\
&=&
2\pi i um^{1-d}
\int_0^\infty
(uy)^{1/m-(\ell+d)/m}
\\
&\times&
\sum_{n=0}^\infty
\frac{(-1)^n}{n!\Gamma(n+1-(\ell+d))}
\big(m(uy)^{1/m}\big)^{2n-(\ell+d)}
\psi(y)dy.
\end{eqnarray*}
By the series definition of the Bessel-function
in (3.1), we have
$$
m\mathfrak{I}_d
=
2\pi i um^{1-d}
\int_0^\infty
(uy)^{1/m-(\ell+d)/m}
J_{-(\ell+d)}\big(2m(uy)^{1/m}\big)
\psi(y)dy.
$$
Applying the asymptotic expansion in (4.2) with
$z=2m(uy)^{1/m}=2\pi m(xy)^{1/m}$, $k=\ell+d$, we get
\begin{eqnarray*}
m\mathfrak{I}_d
&=&
i\pi^{1+\ell-d} m^{1/2-d}x
\sum_{0\leq j<2L}
\frac{\Gamma\big(\ell+d+j+\frac12\big)}
{j!\Gamma\big(\ell+d-j+\frac12\big)(4\pi m)^j}
\\
&\times&
\int_0^\infty
(xy)^{1/(2m)-(\ell+d+j)/m}
\\
&\times&
\Big\{i^{j+\ell+d-1/2}e\big(m(xy)^{1/m}\big)
+(-i)^{\ell+d+j-1/2}e\big(-m(xy)^{1/m}\big)
\Big\}
\psi(y)dy
\\
&+&
O\big((xX)^{-2L/m+1/(2m)-d/m+1/2}\big).
\end{eqnarray*}
Let $L=[(r-d)/2]+1$. Since $L\geq(r-d)/2+1/2$, the error
term above is $O((xX)^{-r/m+1/2})$, and this finishes the estimate for $m\mathfrak{I}_d.$

To estimate $\mathfrak{I}_r^*$, we move the integral
contour from Re$(s)=\sigma_1(r)$ in (2.17)
to Re $(s)=\s(r)=1/4+r/(2m)-\varepsilon$ in (2.5). Then Re $(ms-(\ell+r-1/2))$
moves from $m\varepsilon$ to $-(\ell+r-1)/2-m\varepsilon$.
Picking up residues in this region at
$ms-(\ell+d-1/2)=-n$ for $0\leq n\leq(\ell+r-1)/2$, we get
\begin{eqnarray}
\mathfrak{I}_r^*
&=&
2\pi i u
\sum_{n=0}^{[(\ell+r-1)/2]}
\frac{(-1)^nm^{2(n-r)-\ell+1}}
{\big(\ell+r-n-\frac12\big)n!\Gamma\big(n-\ell-r+1\big)}
\nonumber
\\
&\times&
\int_0^\infty
(uy)^{1/m+2(n-\ell-r)/m}
\psi(y)dy
\nonumber
\\
&+&
\int_{\text{Re}\, s=\s(r)}
\frac{\Gamma\big(ms-\big(\ell+r-\frac12\big)\big)}
{s\Gamma\big(-ms+\frac12\big)}
m^{-2ms+m/2}u^{-2s+1}
\Tilde{\psi}(-2s+1)ds.
\end{eqnarray}
By Stirling's formula, for $\text{Re}\, s=\s(r)$ one has
\begin{eqnarray*}
\frac{\Gamma\big(ms-\big(\ell+r-\frac12\big)\big)}
{s\Gamma\big(-ms+\frac12\big)}
\ll
|s|^{2m\s(r)-\ell-r-1}
\ll
|s|^{-1-\varepsilon}.
\end{eqnarray*}
By (2.18) the last integral in (4.4) is
$O((xX)^{-{r}/m+1/2+2\varepsilon}).$
The first quantity on the right of (4.4) is
\begin{eqnarray*}
\ll (uX)^{1+\frac1{m}+\frac{2[(\ell+r-1)/2]-2\ell-2r}{m}}\ll (xX)^{-r/m+1/2}.
\end{eqnarray*}
Back to (4.3) and collecting our results on $m\mathfrak I_d$,
$\mathfrak I_r^*$ and their error terms, we get
\begin{eqnarray*}
\mathcal {H}_{1}
&=&
-\frac{\sqrt{m}x}{\Gamma\big(\ell+\frac12\big)}
\sum_{b=0}^{r-1}
(-1)^{b}\Gamma\Big(\ell+b+\frac12\Big)
(\pi m)^{-b-1}
\\
&\times&
\sum_{j=0}^{2[(r-b-1)/2]+1}
\frac{\Gamma\big(\ell+b+1+j+\frac12\big)}
{j!\Gamma\big(\ell+b+1-j+\frac12\big)(4\pi m)^j}
\int_0^\infty
(xy)^{1/(2m)-(\ell+b+1+j)/m}
\\
&\times&
\Big\{i^{\ell+b+j+1/2}e\big(m(xy)^{1/m}\big)
+(-i)^{\ell+b+j+1/2}e\big(-m(xy)^{1/m}\big)\Big\}
\psi(y)dy
\\
&+&
O\big((xX)^{-r/m+1/2+\varepsilon}\big).
\end{eqnarray*}
We write $t=b+j+1$ and collect like terms with coefficients
\begin{eqnarray*}
b_t^{(1)}
&=&
-\frac{4\sqrt{m}}{(4\pi m)^t}\cdot
\frac{\Gamma\big(\ell+t+\frac12\big)}
{\Gamma\big(\ell+\frac12\big)}
\sum_{0\leq b<\min(r,t)}
\frac{(-4)^{b}\Gamma\big(\ell+b+\frac12\big)}
{\Gamma(t-b)\Gamma\big(\ell+2b-t+\frac52\big)},
\\
t
&=&
1, 2, \ldots,r+1.
\end{eqnarray*}
Then we obtain
\begin{eqnarray*}
\mathcal {H}_{1}
&=&
x\sum_{t=1}^{r+1}
b_t^{(1)}
\int_0^\infty
(xy)^{1/(2m)-(t+\ell)/m}
\\
&\times&
\Big\{i^{\ell+t-1/2}e\big(m(xy)^{1/m}\big)
+(-i)^{\ell+t-1/2}e\big(-m(xy)^{1/m}\big)\Big\}
\psi(y)dy
\\
&+&
O\big((xX)^{-r/m+1/2+\varepsilon}\big).
\end{eqnarray*}

One can estimates $\mathcal {H}_{j}$ ($j\geq 2$) in a similar way, and  finally obtain
\begin{eqnarray*}
\mathcal {H}_{j}
&=&
x\sum_{t=j}^{r+1}
b_t^{(j)}
\int_0^\infty
(xy)^{1/(2m)-(t+\ell)/m}
\\
&\times&
\Big\{i^{\ell+t-1/2}e\big(m(xy)^{1/m}\big)
+(-i)^{\ell+t-1/2}e\big(-m(xy)^{1/m}\big)\Big\}
\psi(y)dy
\\
&+&
O\big((xX)^{-r/m+1/2+\varepsilon}\big).
\end{eqnarray*}
Substituting these estimates into (2.19) and putting the term with $t=r+1$ into the error term we finish the
proof of Theorem 1.1 for even $m$ by letting $c_k=\sum_{j=0}^kb_k^{(j)}$ where $c_0=b_0^{(0)}=-1/\sqrt{m}$.

\section{Proof of Theorem 1.2}
\setcounter{equation}{0}

By (1.9) with $\psi(x)=\phi(x/X)e(\alpha x^{\beta})$, we get
\begin{eqnarray}
\sum_{n\neq0}
A_f(n,1,\ldots,1)e(\alpha |n|^{\beta})\phi\Big(\frac{|n|}{X}\Big)
=\sum_{n\neq0}
\frac{A_f(1,\ldots, 1, n)}{|n|}{\Phi}\big(|n|\big),
\end{eqnarray}
where by Theorem 1.1, for $r>m/2$ one has
\begin{eqnarray*}
{\Phi}(x)
&=&
x\sum_{k=0}^{r}
\tilde{c}_k\sum_{\pm}  (\pm i)^{k+(m-1)/2}
\int_0^\infty (xy)^{1/(2m)-1/2-k/m}
\phi\Big(\frac{y}{X}\Big)
e\big(\alpha y^{\beta}\pm m(xy)^{1/m}\big)dy
\\
&+&
O\big((xX)^{-r/m+1/2+\varepsilon}\big).
\end{eqnarray*}
Making change of variable $y=t^mX$ and putting in (5.1), we get
\begin{eqnarray}
&&
\sum_{n\neq0}
A_f(n,1,\cdots,1)e(\alpha |n|^{\beta})\phi\Big(\frac{|n|}{X}\Big)
\nonumber
\\
&=&
m\sum_{k=0}^{r}\tilde{c}_k X^{ 1/(2m)+1/2-k/m}\sum_{n=1}^\infty\frac{A_f(1,\cdots, 1, n)+A_f(1,\cdots, 1, -n)}{
n^{1/2+k/m-1/(2m)}}\nonumber\\
&&\qquad\times \sum_{\pm}  (\pm i)^{k+(m-1)/2}I_k(n;\pm)
\nonumber
\\
&+&
O\Big(X^{-r/m+1/2+\varepsilon}\sum_{n=1}^\infty\frac{|A_f(1,\cdots,1,n)|+|A_f(1,\cdots,1,-n)|}{n^{r/m+1/2-\varepsilon}}\Big),
\end{eqnarray}
where
\begin{eqnarray}
I_k(n; \pm)=\int_0^\infty t^{m/2-k-1/2}
\phi(t^m)e\big(\alpha X^\beta t^{m\beta}\pm (nX)^{1/m}mt\big)dt.
\end{eqnarray}
By Rankin-Selberg method for $GL(n)\times GL(n)$ convolution \cite{Gdf} (Remark\ 12.1.8), one has
\begin{eqnarray}
\sum_{1\leq |n|\leq X}|A_f(1,\cdots,1, n)|^2\ll X.
\end{eqnarray}
Therefore the $O$-term in (5.2) is
$O( X^{-r/m+1/2+\varepsilon})$ for $r>m/2$, where the implied constant depends on $r, m$ and $\ve$.
To estimate the integral in (5.3) we consider integral of the form
$$
\int_0^\infty
h(t)e\big(f(t)\big)dt,
$$
where $h, f\in C^\infty(\mathbb R)$ and $h$ is supported on $[a,b]\subset (0,\infty)$. Suppose $f'(x)\not=0$ for $x\in [a,b]$.
By repeated partial integrating by parts, one obtains for $j\geq 0$ that
\begin{eqnarray}
\int_0^\infty
h(t)e\big(f(t)\big)dt
=\Big(\frac{-1}{2\pi i}\Big)^{j} \int_0^\infty
h_j(t)e\big(f(t)\big)dt,
\end{eqnarray}
where $h_0(t)=h(t)$ and
\begin{eqnarray*}
 h_j(t)=\Big(\frac{h_{j-1}(t)}{f'(t)}\Big)':=\frac{g_j(t)}{(f'(t))^{2j}},\qquad j\geq 1.
\end{eqnarray*}

Let $h(t)=t^{m/2-k-1/2}\phi(t^m)$ and $f(t)=\alpha X^\beta t^{m\beta}+ m(nX)^{1/m}t$ in (5.5) one easily obtains
\begin{eqnarray*}
I_k(n;+)\ll_{m,j} (nX)^{-j/m}, \quad {\rm for}\quad n\geq 1.
\end{eqnarray*}
Set $j=r+1$. Then the contribution of the terms in $\sum_{+}$ to (5.2) is
\begin{eqnarray*}
\ll X^{-1/(2m)+1/2-r/m}\sum_{n=1}^\infty\frac{|A_f(1,\cdots, 1, n)|+|A_f(1,\cdots, 1, -n)|}{n^{r/m+1/2+1/(2m)}}
\ll_{r,m} X^{-1/(2m)+1/2-r/m}.
\end{eqnarray*}
Here we have used (5.4). To estimate the contribution of the terms concerning $\sum_{-}$, we write
\begin{eqnarray}
&&n_0=\frac12\min\{1,2^{\beta-1/m}\}(\a\b X^\beta)^mX^{-1}, \\
&&n_1=2\max\{1,2^{\beta-1/m}\}(\a\b X^\beta)^mX^{-1}.
\end{eqnarray}
Then for $n\not\in (n_0,n_1)$ one has $I_k(n;-)\ll (nX)^{-j/m}$ . Therefore the contribution of the terms with $n\not\in (n_0,n_1)$ in $\sum_{-}$ is $O(X^{-1/(2m)+1/2-r/m}).$
This shows that
\begin{eqnarray}
&&
\sum_{n\neq0}
A_f(n,1,\cdots,1)e(\alpha |n|^{\beta})\phi\Big(\frac{|n|}{X}\Big)
\nonumber
\\
&=&
m\sum_{k=0}^{r}\tilde{c}_k (-i)^{k+(m-1)/2}X^{ 1/(2m)+1/2-k/m}\nonumber
\\
&&\times\sum_{n_0<n<n_1}\frac{A_f(1,\cdots, 1, n)+A_f(1,\cdots, 1, -n)}{
n^{1/2+k/m-1/(2m)}}I_k(n;-)\nonumber\\
&&+O_{r,m,\ve}\Big(X^{-r/m+1/2+\varepsilon}\Big).
\end{eqnarray}

If $2\max\{1,2^{\beta-1/m}\}(\a\beta)^m\leq X^{1-\beta m}$, then $n_1\leq 1$. Hence the main term in (5.8) disappears  and the estimate
\bna
\sum_{n\neq0} A_f(n,1,\cdots,1)e(\alpha |n|^{\beta})\phi\Big(\frac{|n|}{X}\Big)
\ll_{m,\beta,r} X^{-r/m+1/2+\varepsilon}\ll_{m,\beta,M} X^{-M}
\ena
holds for any $M>0$ by taking $r$ sufficiently large in terms of $M$. This proves Theorem 1.2 (i).

If\  $2\max\{1,2^{\beta-1/m}\}(\a\beta)^m> X^{1-\beta m},$ then $n_1> 1$.  We distinguish two cases according to $\beta \not=1/m$ or not.
For $\beta \not=1/m$ we have
$$
(\alpha X^\beta t^{m\beta}-(nX)^{1/m}mt\big)''=\alpha (m\beta) (m\beta-1)X^\beta t^{m\beta-2}\gg_{m,\b} \alpha X^{\beta}.
$$
By the second derivative test one has
$I_k(n;-)\ll_{\beta,m} (\a X^\beta)^{-1/2}$. Thus the main term in (5.8) is
\begin{eqnarray*}
&&\ll_{m,\beta}X^{ 1/(2m)+1/2}(\a X^\beta)^{-1/2}\sum_{n_0<n<n_1}\frac{|A_f(1,\cdots, 1, n)|+|A_f(1,\cdots, 1, -n)|}{n^{1/2-1/(2m)}}\\
&&\ll_{m,\beta}(n_1X)^{ 1/(2m)+1/2}(\a X^\beta)^{-1/2}\ll_{m,\beta}(\a X^\beta)^{m/2}.
\end{eqnarray*}
Choosing $r=[(m+1)/2]$ we get
$$
\sum_{n\neq0} A_f(n,1,\cdots,1)e(\alpha n^{\beta})\phi\Big(\frac{|n|}{X}\Big)\ll_{m,\beta} (\a X^\beta)^{m/2}.
$$

For $\b=1/m$, we use the obvious estimate $I_k(n;-)\ll 1$ in (5.8) to get
\begin{eqnarray*}
\sum_{n\neq0}
A_f(n,1,\cdots,1)e(\alpha |n|^{\beta})\phi\Big(\frac{|n|}{X}\Big)
\ll (n_1 X)^{1/2+1/(2m)}\ll_m (\a X^\b)^{(m+1)/2}.
\end{eqnarray*}
This proves (1.14).

Moreover, when  $\b=1/m$, one has $I=(n_0,n_1)=((\alpha/m)^m/2, 2(\alpha/m)^m)$ with $2(\a/m)^m\geq 1$, and
\bna
I_k(n;-)&=&\int_0^\infty t^{m/2-k-1/2}\phi(t^m)e\big((\alpha- mn^{1/m})X^{1/m}t\big)dt.
\ena
Since $(\a/m)^m>1/2$, there is an unique integer $n_\a\geq 1$  such that
$$
(\a/m)^m=n_\a+\lambda,\quad -1/2<\lambda\leq 1/2.
$$
For $n\in I,\ n\not=n_\a$, one has $\le|n^{1/m}-\a/m\ri|\gg_m \le|n-n_\a\ri|\a^{1-m}$. By repeated partial integrating by parts we get
\bna
I_k(n;-)\ll_{m,j} \frac{1}{(\le|n-n_\a\ri|\alpha^{1-m}X^{1/m})^j},\quad j\geq 0.
\ena
Putting in (5.8) and applying (5.4),  the main terms except the term with $n=n_{\alpha}$ produce the contribution which is, for $j\geq 1$,
\begin{eqnarray}
&\ll_m&X^{ 1/(2m)+1/2}(\a^{m-1}X^{-1/m})^j\sum_{n_0<n<n_1\atop n\not=n_{\alpha}}\frac{|A_f(1,\cdots, 1, n)|+|A_f(1,\cdots, 1, -n)|}{
n^{1/2-1/(2m)}}\cdot\frac1{\le|n-n_\alpha \ri|^j}\nonumber\\
&\ll_m&X^{ 1/(2m)+1/2}(\alpha^{m-1} X^{-1/m})^{j}n_1^{1/(2m)}\ll_mX^{ 1/(2m)+1/2}(\alpha^{m-1} X^{-1/m})^{j}\alpha^{1/2}.
\end{eqnarray}
Here we have used (5.4).
Let $0<\varepsilon<1/m$. If $X>\alpha^{m (m-1)/(1-m\varepsilon)}$, one has  $\alpha^{m-1}X^{-1/m}< X^{-\varepsilon}$.
Thus  the last expression in (5.9) is $\ll  X^{-r/m+1/2+\varepsilon}$ by taking $j$ sufficiently large in terms of $r$. This proves
\begin{eqnarray}
&&\sum_{n\neq0}
A_f(n,1,\cdots,1)e(\alpha |n|^{\beta})\phi\Big(\frac{|n|}{X}\Big)
\nonumber\\
&&=m\left\{A_f(1,\cdots, 1, n_\alpha)+A_f(1,\cdots, 1, -n_\alpha)\right\}\sum_{k=0}^{r}\rho_{+}(k,m,\a, X)X^{ 1/(2m)+1/2-k/m}\nonumber\\
&&+O_{r,m, \ve}\big( X^{-r/m+1/2+\varepsilon}\big)
\end{eqnarray}
with
$$
\rho_{+}(k,m,\a, X)=\tilde{c}_k (-i)^{k+(m-1)/2}\frac{I_k(n;-)}{n_{\alpha}^{1/2+k/m-1/(2m)}}.
$$

In particular, suppose $(\alpha/m)^m=q$ is an integer, that is $\a=mq^{1/m}$. Then $n_{\a}=q$ and
\begin{eqnarray*}
I_k(n;-)=\int_0^\infty t^{m/2-k-1/2}\phi(t^m)=\frac{1}{m}
\int_0^\infty x^{1/(2m)-1/2-k/m}\phi(x)dx.
\end{eqnarray*}

This finishes the proof of Theorem 1.2 and  Corollary 1.1 with the exponential function being $e(\a |n|^{\b})$.  Proof for the case of $e(-\a |n|^{\b})$ is analogous.

\section{Proof of Theorems 1.3 and 1.4}
\setcounter{equation}{0}

Let $\Delta>1$  and $\phi:\mathbb R^+\to [0,1]$ be a $C^\infty$ function supported on $[1-\Delta^{-1}, 2+\Delta^{-1}]$ such that $\phi(x)\equiv 1$ for $x\in [1,2]$ and
 satisfies
\begin{eqnarray}
\phi^{(j)}(x)\ll \Delta^{j},\qquad {\rm for\ any\ integer}\quad j\geq 0.
\end{eqnarray}
Then by (5.4) and Cauchy's inequality we get
\begin{eqnarray}
&&\sum_{X<|n|\leq 2X}A_f(n,1,\cdots,1)e(\a |n|^\b)\nonumber\\
&&=\sum_{n\neq0}A_f(n,1,\cdots,1)e(\a |n|^\b)\phi\Big(\frac{|n|}{X}\Big)+
O\big(X\Delta^{-1/2}\big).
\end{eqnarray}
By (5.2) and applying (5.4) again we have
\begin{eqnarray}
&&
\sum_{n\neq0}
A_f(n,1,\cdots,1)e(\alpha |n|^{\beta})\phi\Big(\frac{|n|}{X}\Big)
\nonumber
\\
&=&
m\sum_{k=0}^{r}\tilde{c}_k X^{ 1/(2m)+1/2-k/m}\sum_{n>0}\frac{A_f(1,\cdots, 1, n)+A_f(1,\cdots, 1, -n)}{
n^{1/2+k/m-1/(2m)}}\nonumber\\
&&\times\sum_{\pm}  (\pm i)^{k+(m-1)/2}I_k(n; \pm)+O_m(1),
\end{eqnarray}
where $r=[(m+1)/2]$ and $I_k(n; \pm)$ is defined as in (5.3).
By (5.5) and (6.1), for $j\geq 1$  we have
\begin{eqnarray*}
&&I_k(n; +)\ll_{m,k,j} (nX)^{-j/m}\Delta^{j-1},\quad {\rm for \ }\quad n\geq 1.
\end{eqnarray*}
Set $j=r$ for $n>H=\D^mX^{-1}$ and $j=1$ for $n\leq H$. The contribution of $I(n; +)$ to (6.3) is
\begin{eqnarray*}
&\ll& m X^{1/2-1/(2m)}\sum_{n\leq H}\frac{|A_f(1,\cdots, 1, n)|+|A_f(1,\cdots, 1, -n)|}{n^{1/2+1/(2m)}}\\
&&+mX^{ 1/(2m)+1/2-r/m}\D^{r-1}\sum_{n>H}\frac{|A_f(1,\cdots, 1, n)|+|A_f(1,\cdots, 1, -n)|}{
n^{1/2+r/m-1/(2m)}}
\nonumber\\
&\ll& m (XH)^{ 1/2-1/(2m)}+m(XH)^{ 1/(2m)+1/2-r/m}\D^{r-1}\\
&\ll&m\D^{(m-1)/2}.
\end{eqnarray*}

Next, let $n_1$ be given by (5.7). Then  for $j\geq 1$,
\begin{eqnarray}
&&I_k(n; -)\ll_{m,k,j} (nX)^{-j/m}\Delta^{j-1},\quad {\rm for \ }\quad n\geq n_1.
\end{eqnarray}
Applying (6.4) with $j=1$ for $n_1\leq n\leq H=\D^mX^{-1}$ and $j=r$ for $n>H$, then
the contribution of $I_k(n; -)$ with $n\geq n_1$ to (6.3) is $O(\D^{(m-1)/2})$. This together with (6.2) shows that
\begin{eqnarray}
&&\sum_{X<|n|\leq 2X}A_f(n,1,\cdots,1)e(\a |n|^\b)\nonumber\\
&=&
m\sum_{k=0}^{r}\tilde{c}_k X^{ 1/(2m)+1/2-k/m}\sum_{1\leq n<n_1}\frac{A_f(1,\cdots, 1, n)+A_f(1,\cdots, 1, -n)}{
n^{1/2+k/m-1/(2m)}}\nonumber\\
&&\times\sum_{\pm}  (\pm i)^{k+(m-1)/2}I_k(n; \pm)+O(\D^{(m-1)/2})+O\big(X\Delta^{-1/2}\big).
\end{eqnarray}

Suppose that the parameters $\a, \b, X$ satisfy (1.11). Then $n_1<1$ and the main term above disappears. Setting $\D=X^{2/m}$ we get
\begin{eqnarray*}
\sum_{X<|n|\leq 2X}A_f(n,1,\cdots,1)e(\a |n|^\b)
\ll_m \D^{(m-1)/2}+\D X^{-1/2}\ll_m X^{1-1/m}.
\end{eqnarray*}

Suppose that the parameters $\a, \b, X$ satisfy (1.13), then $n_1>1$. Now we have
\begin{eqnarray}
I_k(n; -)&=&\int_{1}^{2^{1/m}}t^{m/2-k-1/2}e\big(\alpha X^\beta t^{m\beta}-(nX)^{1/m}mt\big)dt+O(\D^{-1}).
\end{eqnarray}
For $\b\not=1/m$, the second derivative test shows that
\begin{eqnarray*}
I_k(n; -)\ll_{m,\b,k} (\a X^\b)^{-1/2}+\D^{-1}.
\end{eqnarray*}
Thus the main term of (6.5) is
\begin{eqnarray*}
&\ll&X^{1/2+1/(2m)}\le(\D^{-1}+(\a X^\b)^{-1/2}\ri)\sum_{1\leq n\leq n_1}\frac{|A_f(1,\cdots, 1, n)|+|A_f(1,\cdots, 1, -n)|}{
n^{1/2-1/(2m)}}\\
&\ll& (Xn_1)^{1/(2m)+1/2}\le(\D^{-1}+(\a X^\beta)^{-1/2}\ri)\\
&\ll& (\a X^\b)^{(m+1)/2}\D^{-1}+(\a X^\b)^{m/2}.
\end{eqnarray*}
Setting $\D=\max\{\a X^{\b}, X^{2/m}\}$,  we yield
\begin{eqnarray*}
\sum_{X<|n|\leq 2X}A_f(n,1,\cdots,1)e(\a |n|^\b)\ll  (\a X^\b)^{m/2}+X^{1-1/m}.
\end{eqnarray*}

For $\b=1/m$, we let $n_\a\geq 1$ be the integer such that
$
(\a/m)^m-n_\a\in (-1/2, 1/2].
$
By (6.6) for $n\not=n_{\a}$,
\bna
I_k(n;-)\ll_{m} \frac{1}{|n-n_\a|\alpha^{1-m}X^{1/m}}+O(\Delta^{-1}).
\ena
Choosing $\D=\max\{\a X^{\b}, X^{2/m}\}$, the  contribution of the terms with $n\not=n_\a$ in (6.5) is
\begin{eqnarray*}
&\ll&X^{1/2-1/(2m)}\a^{(m-1)}\sum_{1\leq n\leq n_1\atop n\not=n_{\alpha}}\frac{|A_f(1,\cdots, 1, n)|+|A_f(1,\cdots, 1, -n)|}{
n^{1/2-1/(2m)}}\cdot\frac1{|n-n_\alpha |}\\
&&+X^{ 1/(2m)+1/2}\Delta^{-1}\sum_{1\leq n\leq n_1\atop n\not=n_{\alpha}}\frac{|A_f(1,\cdots, 1, n)|+|A_f(1,\cdots, 1, -n)|}{
n^{1/2-1/(2m)}}\\
&\ll&\a^{m-1/2}X^{1/2-1/(2m)}+(n_1X)^{1/2+1/(2m)}\Delta^{-1}\nonumber\\
&\ll& \a^{m-1/2}X^{1/2-1/(2m)}.
\end{eqnarray*}
Thus we get
\begin{eqnarray}
&&\sum_{X<|n|\leq 2X}A_f(n,1,\cdots,1)e(\a |n|^\b)\nonumber
\\
&&=m\tilde{c}_0 X^{ 1/(2m)+1/2}\frac{A_f(1,\cdots, 1, n_\a)+A_f(1,\cdots, 1, -n_\a)}{n_\a^{1/2-1/(2m)}}(-i)^{(m-1)/2}I_0(n_\a; -)
\nonumber
\\
&&+ O_m(\a^{m-1/2}X^{1/2-1/(2m)})+O_m(X^{1-1/m}).
\end{eqnarray}
By (5.4) the above fraction is $O(\a^{1/2})$.
Therefore the main term in (6.7) is $O(\a^{1/2}X^{ 1/(2m)+1/2})$. This proves (1.18) and hence finishes the proof of Theorem 1.3.

Now we assume the following
$$
A_f(1,\cdots, 1, n)\ll n^{\th},\quad \rm{for\ some}\quad 0<\th<\frac2{m-1}.
$$
Then the error term $O(X\D^{-1/2})$  in (6.2) and (6.5) becomes $O(X^{1+\th}\D^{-1})$.
Suppose that the parameters $\a, \b, X$ satisfy (1.11), then we can take $\D=X^{2(1+\th)/(m+1)}$ and  get
\begin{eqnarray*}
\sum_{X<|n|\leq 2X}A_f(n,1,\cdots,1)e(\a |n|^\b)
\ll_m X^{(m-1)(1+\th)/(m+1)}.
\end{eqnarray*}
Suppose the parameters $\a, \b, X$ satisfy (1.13), then we can set $\D=\max\{\a X^{\b}, X^{2(1+\th)/(m+1)}\}$ to yield, for $\b\not=1/m$,
\begin{eqnarray*}
\sum_{X<|n|\leq 2X}A_f(n,1,\cdots,1)e(\a |n|^\b)\ll  (\a X^\b)^{m/2}+X^{(m-1)(1+\th)/(m+1)},
\end{eqnarray*}
and for $\b=1/m$,
\begin{eqnarray}
&&\sum_{X<|n|\leq 2X}A_f(n,1,\cdots,1)e(\a |n|^\b)\nonumber
\\
&&=m\tilde{c}_0 X^{ 1/(2m)+1/2}\frac{A_f(1,\cdots, 1, n_\a)+A_f(1,\cdots, 1, -n_\a)}{n_\a^{1/2-1/(2m)}}(-i)^{(m-1)/2}I_0(n_\a; -)
\nonumber
\\
&&+ O_m(\a^{m-1/2}X^{1/2-1/(2m)})+O_m(X^{(m-1)(1+\th)/(m+1)}).
\end{eqnarray}
Note that $I_0(n_\a; -)=I(m,\a, X)+O(\D^{-1}),$
where
$$
I(m,\a, X)=\int_{1}^{2^{1/m}}t^{m/2-1/2}e\big((\alpha-mn_\a^{1/m} )X^{1/m}t\big)dt.
$$
Thus the main term in (6.8) can be rewritten as
$$
m\tilde{c}_0 X^{ 1/(2m)+1/2}\frac{A_f(1,\cdots, 1, n_\a)+A_f(1,\cdots, 1, -n_\a)}{n_\a^{1/2-1/(2m)}}(-i)^{(m-1)/2}I(m,\a,X).
$$
This finishes the proof of Theorem 1.4.

\section{Conclusion and discussion}
\setcounter{equation}{0}

Two important features of smoothly weighted sums of Fourier coefficients of a Maass form for $GL_m(\mathbb Z)$ against $e(\alpha n^{\beta})$
have been discovered. They are rapid decay and resonance for various combinations of $\alpha$ and $\beta$.
These features capture the vibration behavior of the Fourier coefficients of a Maass form.
It is interesting to see whether these two features can be used to characterize Fourier coefficients of Maass forms.
On the other hand, when $\beta$ is large, our methods failed to derive a non-trivial bound for the smoothly weighted sum.
These are subjects of our subsequent research.

%{\bf Acknowledgment.} The first author was supported by the National Natural Science Foundation of China (Grant No. 10971119).

\end{document}